\documentclass[10pt,leqno]{amsart}
\usepackage{srcltx}
\usepackage{latexsym,mathrsfs}

\usepackage{amsmath,amssymb}
\usepackage{amsxtra}
\usepackage{amsfonts}
\usepackage{amsthm}
\usepackage{newlfont}
\newtheorem{theorem}{Theorem}[section]
\newtheorem{cor}[theorem]{Corollary}
\newtheorem{lemma}[theorem]{Lemma}

\def\R{{\mathbb{R}}}

\newcommand{\ad}{\mathop{\mathrm{ad}}\nolimits}

\newcommand{\Erf}{\mathop{\mathrm{Erf}}\nolimits}

\newcommand\Pp{{\mathtt P}}

\numberwithin{equation}{section} 

\title[]{Orthogonal matrix polynomials satisfying differential equations with recurrence coefficients having non-scalar limits}

\author{Jorge Borrego}
\address{J. Borrego \\
Departamento de Matem\'aticas\\
Universidad Carlos III de Madrid \\
Avda. de la Universidad 30, 28911 Legan\'es, Madrid, Spain.}
\email{jborrego@math.uc3m.es}
\author{Mirta Castro}
\address{M. Castro \\
Departamento de Matem\'atica Aplicada II\\
Universidad de Sevilla \\
E.P.S., c/Virgen de \'Africa 7, 41011 Sevilla, Spain.}
\email{mirta@us.es }
\author{Antonio J. Dur\'an}
\address{A. J. Dur\'an \\
Departamento de An\'{a}lisis Matem\'{a}tico \\
Universidad de Sevilla \\
Apdo (P. O. BOX) 1160\\
41080 Sevilla. Spain.}
\email{duran@us.es }
\thanks{The
     work of the authors is  partially supported by MCI ref.
     MTM2009-12740-C03-01 (first author) and MCI ref. MTM2009-12740-C03-02, (\textit {Junta de
Andaluc\'{\i}a}) ref. FQM-262, FQM-4643 (second and third
authors).}



\date{}

\begin{document}
\maketitle

\begin{abstract}
We introduce a  family of weight matrices $W$  of the form
$T(t)T^*(t)$, $T(t)=e^{\mathscr{A}t}e^{\mathscr{D}t^2}$, where
$\mathscr{A}$ is certain nilpotent matrix and $\mathscr{D}$ is a diagonal matrix
with negative real entries. The weight matrices $W$ have arbitrary size
$N\times N$ and depend on $N$ parameters.

The orthogonal polynomials with respect to this family of weight matrices
satisfy a second order differential equation with
differential coefficients that are matrix polynomials $F_2$, $F_1$ and $F_0$ (independent of $n$)
of degrees not bigger than $2$, $1$ and $0$ respectively.

For size $2\times 2$, we find an explicit
expression for a sequence of orthonormal polynomials with
respect to $W$. In particular, we show that one of the recurrence coefficients
for this sequence of orthonormal polynomials does not asymptotically behave
as a scalar multiple of the identity, as it happens in the examples
studied up to now in the literature.

\end{abstract}

\section{Introduction}

In the last few years a large class of families of $N\times N$
weight matrices $W$ having  symmetric second order differential operators of
the form
\begin{align}\label{difope}
\left( \frac{d}{dt}\right)^2F_2(t)+\left( \frac{d}{dt}\right)^1F_1(t)+F_0(t)
\end{align}
has been introduced (see \cite{DG1}, \cite{DG5}, \cite{D2}, \cite{D3},
\cite{DdI1}, \cite{DdI2}, \cite{G}, \cite{GPT2}, \cite{GPT4}, \cite{CMV},  \cite{PT}). The
differential coefficients
   $F_2$, $F_1$ and $F_0$ are matrix polynomials (which do not depend on $n$)
of degrees less
than or equal to $2$, $1$ and $0$, respectively. As usual, the symmetry of an operator $D$ with respect to the
weight matrix $W$ is defined by $\int D(P)dWQ^*=\int PdW(D(Q))^*$, for any matrix polynomials $P,Q$.

A sequence $(P_n)_n$ of orthogonal  polynomials with respect to a weight matrix $W$ is
 a sequence of matrix polynomials satisfying that $P_n$, $n\ge 0$, is a matrix polynomial of degree $n$
with non singular leading coefficient, and $\int P_ndWP_m^*=\Delta_n\delta _{n,m}$, where
$\Delta _n$, $n\ge 0$, is a positive definite matrix.
When $\Delta_n=I$, we say that the polynomials $(P_n)_n$ are orthonormal (we denote by $I$
the identity matrix).

Just as in the scalar case, any sequence $(P_n)_n$ of orthonormal polynomials with respect to a weight matrix
satisfies a three term recurrence relation
\begin{equation}\label{ttrr}
tP_n(t)=A_{n+1}P_{n+1}(t)+B_nP_n(t)+A_n^*P_{n-1}(t), \quad n\ge 0,
\end{equation}
where $A_n$, $n\ge 1$, are nonsingular and $B_n$, $n\ge 0$, Hermitian.

If $(P_n)_n$ is a sequence of
orthonormal matrix polynomials with respect to $W$, the
symmetry of the second order
differential operator (\ref{difope}) is equivalent to the
second order differential equation
\begin{align}\label{difequ}
P''_n(t)F_2(t)+P'_n(t)F_1(t)+P_n(t)F_0(t)=\Lambda _n P_n(t),
\end{align}
where $\Lambda _n$ are
Hermitian matrices (see Lemma 4 of \cite{D1}).

The theory of matrix valued orthogonal
polynomials was started by M. G. Krein in 1949 \cite{K1,K2} (see also \cite{Be} or \cite{Atk}), but more than 50 years have been
necessary to see the first examples  of orthogonal matrix polynomials satisfying that kind of differential
equations (see \cite{DG1,G,GPT2}). These examples will likely play, in the case of matrix orthogonality, the role of the classical families
of Hermite, Laguerre and Jacobi in the case of scalar
orthogonality.

As their scalar relatives, these families of orthogonal matrix polynomials
also satisfy many formal properties, relationships and structural formulas (see \cite{DG3,DL,D4}).
 These structural formulas have been very useful to compute explicitly the orthonormal polynomials related to several of these families, and, in particular,
their re\-cu\-rren\-ce coefficients (\ref{ttrr}). The recurrence coefficients of these examples
asymptotically behave as scalar multiples of the identity. More precisely, either $\lim _n A_n=aI$ and
$\lim _n B_n=bI$ (\cite{DL,GdI}), or $\lim _n a_nA_n=aI$, $a\not =0$, and
$\lim _n b_nB_n=bI$, for certain divergent sequences $(a_n)_n$, $(b_n)_n$ of positive real numbers (\cite{DG3,DL,D4}).

The purpose of this paper is to introduce a new family of weight matrices ha\-ving orthonormal polynomials satisfying
 second order differential equations and whose recurrence coefficients do not asymptotically behave as scalar multiples of the identity.

Our example is the family of weight matrices of arbitrary size $N\times N$ constructed from the $N-1$ non-null complex parameters
$a_1,\cdots a_{N-1}$ and the positive real parameter $b$ ($b\neq 1$) as follows.
Consider the nilpotent matrix $A$ and the diagonal matrices $\mathscr{J}$ and $\Psi$ defined by
\begin{align}\label{defAJ}
A&=\begin{pmatrix}
        0 & a_1   & 0   & \cdots  &  0    \\
        0 &  0    & a_2 &  \cdots &  0    \\
       \vdots &\vdots &\vdots &\ddots &\vdots                       \\
        0 & 0     &   0 &  \cdots& a_{N-1}   \\
        0 & 0     & 0   & \cdots      & 0
      \end{pmatrix}  \quad \mathscr{J}=\begin{pmatrix}
        0 & 0 &0    & \cdots  &  0    \\
        0 &  1&0     &  \cdots &  0    \\
         0 &  0&2     &  \cdots &  0    \\
       \vdots  &\vdots&\vdots &\ddots &\vdots                       \\
        0      & 0 &0  &\cdots     & N-1
      \end{pmatrix},\\ \label{psi2}
      \Psi&=I+\frac{b-1}{N-1}\mathscr{J}.
\end{align}
Let the diagonal matrix $\mathscr{D}$ and the upper triangular nilpotent matrix $\mathscr{A}$
be defined by
\begin{align}\label{D}
\mathscr{D}=-\frac{b}{2}\Psi ^{-1},\\
\label{Acurs}
\mathscr{A}=\displaystyle \sum_{j=0}^{\left[\frac{N}{2}\right]-1}\alpha_jA^{2j+1},
\end{align}
where
$$
\alpha_j=\frac{(1-b)^j(2j+1)^{j-1}}{(4b)^j(N-1)^jj!},\quad j\ge 0.
$$
Our weight matrix $W$ is then defined by
\begin{equation}\label{defW}
W(t)=T(t)T^*(t), \quad T(t)=e^{\mathscr{A}t}e^{\mathscr{D}t^2}.
\end{equation}
Since $\mathscr{A}$ is nilpotent of order $N$, $e^{\mathscr{A}t}$ is always a polynomial of degree $N-1$.

For $b=1$ (considering $\alpha_0=1$), we recover the example 5.1 of \cite{DG1}.

For the benefit of the reader, we display here our weight matrix for size $2\times 2$. For $N=2$ we have
$$
\mathscr{A}=A=
\left(
\begin{array}{cc}
 0 & a  \\
 0 & 0
\end{array}
\right),\quad  \mathscr{D}=\left(
\begin{array}{cccc}
  -\frac{b}{2}   & 0               \\
 0                 & -\frac{1}{2}
\end{array}
\right)$$
and then
\begin{equation}\label{defW2}
W=\left(
\begin{array}{cc}
 |a|^2 t^2e^{-t^2} +e^{-bt^2}  & a  te^{-t^2} \\
 \overline{a} t e^{-t^2}        & e^{-t^2}
\end{array}
\right),\quad T=\left(
\begin{array}{cc}
 e^{-bt^2/2} & ate^{-t^2/2}  \\
 0 & e^{-t^2/2}
\end{array}
\right) .
\end{equation}
The content of this paper is the following. In Section 3, we prove that our weight matrix $W$ always has
a symmetric second order differential operator like (\ref{difope}):

\begin{theorem}\label{th1}
The second order differential operator  (\ref{difope}) with differential coeffi\-cients $F_2$,
$F_1$ and $F_0$ given by
\begin{align}\label{F2}
F_2(t)&=\Psi +\frac{b-1}{N-1}\left[ \mathscr{A},\mathscr{J}\right]t,\\
\label{F1}
F_1(t)&=2\mathscr{A}\Psi +2\left( -bI+\frac{b-1}{N-1}\mathscr{A}[\mathscr{A},\mathscr{J}]\right)t,\\
\label{F0}
F_0(t)&=2b\mathscr{J}+\mathscr{A}^2\Psi ,
\end{align}
is symmetric with respect to the weight matrix $W$ (\ref{defW}) (as usual $[X,Y]$ denotes the commutator of the
matrices $X,Y$: $[X,Y]=XY-YX$).
\end{theorem}

For $N=2$, these differential coefficients are
$$
F_2(t)=\left(
\begin{array}{cc}
 1 & a(b-1)t  \\
 0 & b
\end{array}
\right), \quad F_1(t)=\left(
\begin{array}{cc}
 -2bt & 2ab  \\
 0 & -2bt\end{array}
\right), \quad F_0(t)=\left(
\begin{array}{cc}
 0 & 0  \\
 0 & 2b
\end{array}
\right).
$$

The rest of the Sections are devoted to study in depth the orthogonal polynomials with respect to our weigh matrix for size $2\times 2$.
The study of the orthogonal polynomials for higher size $N$, $N\ge 3$, remains a challenge.

In Section 4, we construct the following Rodrigues' formula for a sequence of
orthogonal polynomials with respect to the weight matrix $W$ (\ref{defW2}):

\begin{theorem}\label{th2}
Let the function
$P_n$, $n\ge 1$, be defined by
\begin{equation}\label{defRodri}
P_n(t)=(-1)^n\left[e^{- t^2}\left(\begin{array}{cc} \displaystyle{b^{-n} e^{(1-b) t^2}+ \frac{\vert a\vert ^2}{2}(n+2 t^2)}&a t\\
\bar a \left[ 2 t+e^{t^2} \sqrt{\pi}n\Big( \Erf(\sqrt{b} t)-\Erf(t)\Big)\right] &2
\end{array}
\right)\right]^{(n)}W^{-1},
\end{equation}
where  $\Erf$ denotes the error function $\Erf(z)=\displaystyle{\frac{2}{\sqrt{\pi} }}\int_0^{z}e^{-x^2}dx$.
Then $P_n$, $n\ge 1$, is a polynomial of degree $n$ with nonsingular leading coefficient equal to
\begin{equation}\label{gamma}
\Gamma_n=2^n\left(\begin{array}{cc} 1&0\\0&\gamma_n
\end{array}
\right), \quad \gamma_n=2+|a|^2b^{n-\frac{1}{2}}n.
\end{equation}
Moreover, defining $P_0=\left(
\begin{array}{cc}
 1 & 0  \\
 0 & 2
\end{array}
\right)$, $(P_n)_n$ is a sequence of
orthogonal polynomials with respect to $W$.
\end{theorem}

This Rodrigues' formula  allows us to find an explicit expression for the polynomials $(P_n)_n$ in terms
of the Hermite polynomials.

\begin{cor} \label{cor1}
For $n\ge 1$, we have
\begin{equation}\label{PnitHn}
P_n(t)=\left(\begin{array}{cc}\displaystyle b^{-n/2}H_{n}(\sqrt{b}t)&-atb^{-n/2}H_{n}(\sqrt{b}t)+\frac{a}{2}H_{n+1}(t)\\
-2\overline{a}b^{n/2}nH_{n-1}(\sqrt{b}t)&2|a|^2b^{n/2}ntH_{n-1}(\sqrt{b}t)+2H_{n}(t)
\end{array}\right),
\end{equation}
where $H_n$ is the $n$-th Hermite polynomial defined by $H_n(t)=(-1)^n\left( e^{-t^2}\right) ^{(n)}e^{t^2}$.
\end{cor}

In Section 6, using again the Rodrigues' formula (\ref{defRodri}), we find the following three term recurrence relation for a sequence
$(\mathscr{P}_n)_n$ of
orthonormal polynomials with respect to $W$.

\bigskip
\begin{theorem}\label{the3} The sequence of matrix polynomials defined by $\mathscr{P}_{-1}=0$,
$\mathscr{P}_0=\displaystyle{(\pi)^{-\frac{1}{4}}\left(\begin{array}{cc} \displaystyle{\sqrt{\frac{2\sqrt{b}}{\gamma_1}}}&0\\0&1
\end{array}
\right)}$ and
\begin{equation}\label{ttrcon}
t\mathscr{P}_n(t)=A_{n+1}\mathscr{P}_{n+1}(t)+B_n\mathscr{P}_n(t)+A^*_n\mathscr{P}_{n-1}(t),\quad n\geq 0,
\end{equation}
where
\begin{align}\label{ttrrPn}
A_n&=\sqrt{n}\left(\begin{array}{cc} \displaystyle{\sqrt{\frac{\gamma_{n+1}}{2b\gamma_n}}}&0
\\0&\displaystyle{\sqrt{\frac{\gamma_{n-1}}{2\gamma_n}}}
\end{array}
\right),
\\
B_n&=\displaystyle{\frac{b^{\frac{2n-3}{4}}(b+(b-1)n)}{\sqrt{\gamma_n\gamma_{n+1}}}}\left(\begin{array}{cc} 0&a\\\bar{a}&0
\end{array}
\right),
\end{align}
are orthonormal with respect to $W$ (\ref{defW2}) (where the sequence $(\gamma _n)_n$ is defined by (\ref{gamma})).
\end{theorem}

This gives for the recurrence coefficients $(A_n)_n$ the asymptotic behaviour
$$
\underset{n\rightarrow \infty}{\lim}\frac{A_n}{\sqrt{n}} =
\left\{\begin{array}{c}\left(\begin{array}{cc} \displaystyle{\frac{1}{\sqrt{2}}}&0\\0&\displaystyle{\frac{1}{\sqrt{2b}}}
\end{array}
\right)\ \textrm{if}\ b>1\\\left(\begin{array}{cc} \displaystyle{\frac{1}{\sqrt{2b}}}&0\\0&\displaystyle{\frac{1}{\sqrt{2}}}
\end{array}
\right)\ \textrm{for}\ 0<b<1\end{array}\right.
.$$
This limit shows that the recurrence coefficients $(A_n)_n$ do not asymptotically behave as a scalar multiple of the identity,
as it happens in the examples studied up to now in the literature.

\section{Preliminaries}

A weight matrix $W$ is an $N\times N$ matrix of measures supported in the real line satisfying that (1) $W(A)$ is positive semidefinite for any Borel
set $A\in \R$, (2) $W$ has finite moments $\int t^ndW(t)$ of every order, and (3) $\int P(t)dW(t)P^*(t)$ is nonsingular if the leading
coefficient of the matrix polynomial $P$ is nonsingular (all the matrices considered in this paper are
square matrices of size $N\times N$). When all the entries of the matrix $W$
have a smooth density with respect to the Lebesgue measure, we will write $W(t)$ for the matrix whose entries are these densities.
Condition (3) above is necessary and sufficient to guarantee the existence of a sequence $(P_n)_n$ of orthogonal matrix polynomials
with respect to $W$.
Condition (3) above is fulfilled, in particular, when $W(t)$ is positive definite
in an interval of the real line. This is the case of the weight matrix (\ref{defW}) introduced in this paper.

The key concept to study orthogonal matrix polynomials $(P_n)_n$ satisfying second order differential equations of the form
\begin{equation}\label{jj}
P''_n(t)F_2(t)+P'_n(t)F_1(t)+P_n(t)F_0(t)=\Lambda _n P_n(t),
\end{equation}
where the differential coefficients $F_2$, $F_1$ and $F_0$ are matrix polynomials (which do not depend on $n$)
of degrees less
than or equal to $2$, $1$ and $0$, respectively, is that of the symmetry of a differential operator with respect to a weight matrix. Indeed, if we write
$$
D=\left( \frac{d}{dt}\right)^2F_2(t)+\left( \frac{d}{dt}\right)^1F_1(t)+F_0(t),
$$
then $D$ is symmetric with respect to $W$ if and only if the orthonormal matrix polynomials $(P_n)_n$ with
respect to $W$ satisfy (\ref{jj}), where $\Lambda _n$, $n\ge 0$, are Hermitian matrices (see \cite{D1}). If $\Lambda _n$ are not
Hermitian, then the operator $D$ can be decomposed as $D=D_1+iD_2$, where $D_1$ and $D_2$
are  second order differential operators of the form (\ref{difope}) symmetric with respect to $W$  (see \cite{GT}).

The symmetry of a second order differential operator as (\ref{difope}) with respect to a weight matrix can be guaranteed
by a set of differential equations (which is the matrix analogous to the Pearson equation $(f_2w)'=f_1w$ of the scalar case):

\begin{theorem} \label{resum0} For a weight matrix $dW=W(t)dt$, $t\in (a,b)$ ($a$ and $b$ finite or infinite)
and matrix polynomials $F_2, F_1, F_0$ of degrees
not larger than $2, 1$ and $0$, the symmetry
of the second order differential operator
(\ref{difope}) with respect to $W$ follows from the set of equations
\begin{align}\label{ccp}
F_2W &= WF_2^*, \\
\label{eqsym}
 2(F_2W)'-F_1W= WF_1^*,&\quad
(F_2W)''-(F_1W)'+F_0W =WF_0^*,
\end{align}
under the boundary conditions
\begin{align}\label{boundcondit}
\lim _{t\to a^+, b^-}t^nF_2(t)W(t)=0,\quad \lim _{t\to a^+, b^-}t^n[(F_2(t)W(t))'-F_1(t)W(t)]=0, \quad n\ge 0.
\end{align}
\end{theorem}

(See \cite{DG1} or  \cite{GPT2}).

To check that the weight matrix $W$ (\ref{defW}) and the coefficients defined by
(\ref{F2}), (\ref{F1}) and (\ref{F0}) satisfy the differential equations (\ref{eqsym}) we will use the following theorem:

\begin{theorem}\label{resum}
Let $\Omega$  be an open set of the real line.
Let $F_2$, $F_1$, $F$ and $T$ be twice differentiable $N\times N$ matrix functions
on $\Omega$, (with $T(t_0)$ nonsingular for certain $t_0\in \Omega$), and define $\displaystyle W(t)=T(t)T^*(t)$.
Under the assumptions
\begin{align}
\label{hypot1} F_2W &=WF_2^*,  \\
\label{hypot2} T'(t)&=F(t) T(t),\quad \mbox{and} \\
\label{hypot3} F_1&=F_2F+FF_2+F_2',
\end{align}
we have
\begin{enumerate}
\item The weight matrix  $W$ satisfies the
first order differential equation
$$
2(F_2W )'=F_1W +WF_1^*.
$$
\item For a given matrix $F_0$, the second order differential equation
$$
(F_2W )''-(F_1W )'+F_0W =WF_0^*,
$$
holds if and only if the matrix function
\begin{equation}\label{eq4.4}
\chi  = T^{-1}(-FF_2F-F'F_2-FF_2'+F_0)T,
\end{equation}
is Hermitian at each point of $\Omega $.
\end{enumerate}
\end{theorem}

Theorem \ref{resum} is a particular case of Theorem 2.3 of \cite{D3} (the special case
when $F_2$ is a scalar function is Theorem 4.1 of \cite{DG1}).

\bigskip

We will need the following theorem (see \cite{D4}) to find the Rodrigues' formula displayed in Theorem \ref{th2}
for the orthogonal polynomials with respect to the weight matrix
$W$ (\ref{defW2}):

\begin{theorem}\label{throdri}
Let $F_2$, $F_1$ and $F_0$ be matrix polynomials of degrees not
larger than $2$, $1$, and $0$, respectively. Let $W$, $R_n$ be
$N\times N$ matrix functions twice and $n$ times differentiable,
respectively, in an open set $\Omega$ of the real line. Assume
that $W(t)$ is nonsingular for $t\in \Omega$ and that satisfies
the identity (\ref{ccp}), and the differential equations
(\ref{eqsym}). Define the functions $\Pp_n$, $n\ge
1$, by
\begin{align}\label{Rodri}
\Pp_n=R_n^{(n)}W^{-1}.
\end{align}
If for a matrix $\Lambda_n$, the function $R_n$ satisfies
\begin{align}\label{eqXn}
(R_nF_2^*)''-\big( R_n\big[ F_1^*+n(F_2^*)'\big]\big)'+R_n\left[ F_0^*+n(F_1^*)'+\binom{n}{2}(F_2^*)''\right]=\Lambda _n R_n.
\end{align}
then the function $\Pp_n$ satisfies
\begin{align}\label{difequ2}
\Pp''_n(t)F_2(t)+\Pp'_n(t)F_1(t)+\Pp_n(t)F_0(t)=\Lambda _n \Pp_n(t).
\end{align}
\end{theorem}

We will also make use of the following well known formula: for any matrices $X,Y \in \mathbb{C}^{N\times N}$:
\begin{equation}\label{adformula}
e^{X}Y=\left(\underset{n\geq 0}{\sum}\frac{t^n}{n!}\textrm{ad}^{n}_XY \right)e^{X},
\end{equation}
where we use the standard notation
$$
\ad^0_XY=Y,\quad \ad^1_XY=[X,Y]=XY-YX, \quad \ad^2_XY=[X,[X,Y]],
$$
and in general, $\ad^{n+1}_XY=[X,[\ad^n_XY]]$.

\section{Symmetric differential operator}

In this Section we prove Theorem \ref{th1}, that is, the second order differential operator with
coefficients given by (\ref{F2}), (\ref{F1}) and (\ref{F0}) is symmetric with respect to the weight
matrix $W$ (\ref{defW}).

We now list  some technical relations
which we will need in the proof of Theorem \ref{th1} (they will be proved later).

\begin{lemma}\label{sumaciones}
Let the function $F_2$ and the matrices $A$, $\mathscr{A}$, $\Psi$, $\mathscr{D}$ and $\mathscr{J}$  be defined
by (\ref{F2}), (\ref{defAJ}), (\ref{psi2}), (\ref{D}) and (\ref{Acurs}), respectively. Then
\begin{align}
\label{suma2s}
[\mathscr{A},\mathscr{J}]&=\displaystyle \sum_{j=0}^{\left[\frac{N}{2}\right]-1}(2j+1)\alpha_jA^{2j+1}.\\
\label{dfe0}
e^{\mathscr{A}t}\Psi &=F_2e^{\mathscr{A}t}.\\
\label{dfe1}
e^{\mathscr{A}t}\mathscr{D}e^{-\mathscr{A}t}\Psi &=-\frac{b}{2}I-\frac{(b-1)t}{N-1}
e^{\mathscr{A}t}\mathscr{D}e^{-\mathscr{A}t}[\mathscr{A},\mathscr{J}].\\
\label{dfe2}
\Psi e^{\mathscr{A}t}\mathscr{D}e^{-\mathscr{A}t}&=-\frac{b}{2}I-\frac{(b-1)t}{N-1}
[\mathscr{A},\mathscr{J}]e^{\mathscr{A}t}\mathscr{D}e^{-\mathscr{A}t}.\\
\label{suma_2k+1}
\mathscr{A}\left[\mathscr{A},\mathscr{J}\right]&=\frac{2b(N-1)}{1-b}\sum_{j=1}^{\left[\frac{N-1}{2}\right]}\alpha_j\frac{(2j)^j}{(2j+1)^{j-1}}A^{2j}.\\
\label{sumaprinc}
\left[\mathscr{A},\mathscr{J}\right]-\mathscr{A}&=\frac{(1-b)}{2b(N-1)}\mathscr{A}^2\left[\mathscr{A},\mathscr{J}\right] .
\end{align}
\end{lemma}

We are now ready to prove Theorem \ref{th1}.

\begin{proof} (of Theorem \ref{th1})

The symmetry of the second order differential operator  with respect to the weight matrix $W$
will be a consequence of Theorems \ref{resum0} and \ref{resum}. We have to check the boundary conditions
(\ref{boundcondit}) and the three equations
(\ref{ccp}) and (\ref{eqsym}).

To make the proof easier to follow, we proceed in four steps.

\textit{First step: Boundary conditions (\ref{boundcondit}).}
Proof: Since $\mathscr{A}$ is nilpotent, we deduce that
$e^{\mathscr{A}t}$ is a polynomial. The matrix function $T=e^{\mathscr{A}t}e^{\mathscr{D}t^2}$ then decays
exponentially at $\infty$ because the entries of the diagonal matrix $\mathscr{D}$ are negative.
Hence the weight matrix $W=TT^*$ also decays exponentially at $\infty$.
Since $F_2$ and $F_1$ are poly\-no\-mials, it follows straightforwardly that
$t^nF_2W$ and $t^n[(F_2W)'-F_1W]$, $n\ge 0$, have vanishing limits at $ \infty$.

\bigskip

\textit{Second step: $F_2W=WF_2^*$.} Proof: Formula (\ref{dfe0}) of Lemma \ref{sumaciones}  shows that $F_2T=T\Psi$, where $\Psi$ is the diagonal
matrix (real entries) defined by (\ref{psi2}) and $T=e^{\mathscr{A}t}e^{\mathscr{D}t^2}$. Then
$$
F_2TT^*=T\Psi T^*=T(T\Psi)^*=T(F_2T)^*=TT^*F_2^*.
$$
Since $W=TT^*$, we get that $F_2W=WF_2^*$.

\textit{Third step: $2(F_2W)'=F_1W+WF_1^*$}. Proof:
To check the equation $2(F_2W)'=F_1W+WF_1^*$, we
use the first part of Theorem \ref{resum} with $\Omega =\R$.
Hence, we have to prove that $T$ satisfies
$T'(t)=F(t)T(t)$, where $F$ is a solution of the matrix
equation
\begin{equation}\label{F1F}
F_1(t)=F_2(t)F(t)+F(t)F_2(t)+F_2'(t).
\end{equation}

Taking into account that $T=e^{\mathscr{A}t}e^{\mathscr{D}t^2}$, a direct computation gives that
\begin{equation}\label{F}
F(t)=\mathscr{A}+2t e^{\mathscr{A}t}\mathscr{D}e^{-\mathscr{A}t}.
\end{equation}

The definition of $F_2$ (\ref{F2}) and
(\ref{F}) give
\begin{align*}
F_2F+FF_2&=\mathscr{A}\Psi+\Psi\mathscr{A}+\frac{2(b-1)t}{N-1}\mathscr{A}[\mathscr{A},\mathscr{J}]
+2t(\Psi e^{\mathscr{A}t}\mathscr{D}e^{-\mathscr{A}t}+
e^{\mathscr{A}t}\mathscr{D}e^{-\mathscr{A}t}\Psi)\\&\quad \quad + \frac{2(b-1)t^2}{N-1}(e^{\mathscr{A}t}
\mathscr{D}e^{-\mathscr{A}t}[\mathscr{A},\mathscr{J}]+[\mathscr{A},\mathscr{J}]e^{\mathscr{A}t}\mathscr{D}e^{-\mathscr{A}t}).
\end{align*}
Using the definition of $\Psi$ (\ref{psi2}), (\ref{dfe1}) and (\ref{dfe2}) of Lemma \ref{sumaciones} we get
$$
F_2F+FF_2=2\mathscr{A}+\frac{(b-1)}{N-1}(\mathscr{J}\mathscr{A}+\mathscr{A}\mathscr{J})+\frac{2t(b-1)}{N-1}\mathscr{A}[\mathscr{A},\mathscr{J}]-2bIt.
$$
Formula (\ref{F1F}) now follows easily taking into account  the definitions of $F_2$ (\ref{F2}), $F_1$ (\ref{F1}) and $\Psi$ (\ref{psi2}).

\bigskip

\textit{Fourth step: $(F_2W)''-(F_1W)'+F_0W=WF_0^*$}. Proof:

Using 2 of Theorem \ref{resum}, this is equivalent to prove that the matrix
\begin{equation} \label{fchi}
\chi=T^{-1}(-FF_2F-(FF_2)'+F_0)T
\end{equation}
is Hermitian.

We actually will prove  that the matrix function $\chi$ defined by (\ref{fchi}) is diagonal with real entries.

Taking into account that $T(t)=e^{\mathscr{A}t}e^{\mathscr{D}t^2}$, and $\mathscr{D}$ is diagonal, it is enough to prove that the matrix function
\begin{equation}\label{a verificar diagonal}
\xi=e^{-\mathscr{A}t}(-FF_2F-(FF_2)'+F_0)e^{\mathscr{A}t}
\end{equation}
is diagonal.

We first compute $e^{-\mathscr{A}t}(FF_2F)e^{\mathscr{A}t}$.

Taking into account the expression for $F(t)$ in (\ref{F}), that $\mathscr{A}$ and $e^{\mathscr{A}t}$ commute and using
(\ref{dfe0}), one has after straightforward computations
\begin{equation*}
e^{-\mathscr{A}t}(FF_2F)e^{\mathscr{A}t}=\mathscr{A}\Psi\mathscr{A}+
2t(\mathscr{A}\Psi \mathscr{D}+\mathscr{D}\Psi\mathscr{A})+4t^2\mathscr{D}\Psi \mathscr{D}.
\end{equation*}
The definition of $\mathscr{D}$ (\ref{D}) and $\Psi $ (\ref{psi2}) now give
\begin{equation}\label{sumando I}
e^{-\mathscr{A}t}(FF_2F)e^{\mathscr{A}t}=\mathscr{A}^2+\frac{b-1}{N-1}\mathscr{A}\mathscr{J}\mathscr{A}-2bt\mathscr{A}-2t^2b\mathscr{D}.
\end{equation}
We now compute $e^{-\mathscr{A}t}(FF_2)'e^{\mathscr{A}t}$. Using again the definition of $F$ (\ref{F}) and (\ref{dfe0}) of Lemma \ref{sumaciones},
 we have that
\begin{equation*}
(FF_2)'=(\mathscr{A}F_2+2te^{\mathscr{A}t}\mathscr{D}\Psi e^{-\mathscr{A}t})'.
\end{equation*}
The definition of $\mathscr{D}$ (\ref{D}) and $F_2$ (\ref{F2}) give
\begin{align*}
(FF_2)'&=\left( \mathscr{A}\Psi+\frac{(b-1)t}{N-1}\mathscr{A}
[\mathscr{A},\mathscr{J}]-b tI\right) '\\&=\frac{b-1}{N-1}\mathscr{A}
[\mathscr{A},\mathscr{J}]-bI.
\end{align*}
Identity (\ref{suma2s}) of Lemma \ref{sumaciones} shows that $e^{\mathscr{A}t}$ and $[\mathscr{A},\mathscr{J}]$ commute
(they are linear combinations of power of $A$).
One then obtains
\begin{equation}\label{sumando 2}
e^{-\mathscr{A}t}(FF_2)'e^{\mathscr{A}t}=\frac{b-1}{N-1}\mathscr{A}
[\mathscr{A},\mathscr{J}]-bI.
\end{equation}
We finally compute $e^{-\mathscr{A}t}F_0e^{\mathscr{A}t}$. The definition of $F_0$ (\ref{F0}) gives
\begin{align*}
e^{-\mathscr{A}t}F_0e^{\mathscr{A}t}&=e^{-\mathscr{A}t}(2b\mathscr{J}+\mathscr{A}^2\Psi)e^{\mathscr{A}t}\\
&=2be^{-\mathscr{A}t}\mathscr{J}e^{\mathscr{A}t}+\mathscr{A}^2e^{-\mathscr{A}t}\Psi e^{\mathscr{A}t}.
\end{align*}
Again (\ref{suma2s}) of Lemma \ref{sumaciones} shows that $\mathscr{A}$ and $[\mathscr{A},\mathscr{J}]$ commute. Hence
$\ad ^n_\mathscr{A}\mathscr{J}=0$, $n\ge 2$.
Using this fact and
(\ref{adformula}) one obtains
\begin{equation*}
e^{-\mathscr{A}t}F_0e^{\mathscr{A}t}=
2b(\mathscr{J}-[\mathscr{A},\mathscr{J}]t)+\mathscr{A}^2(\Psi-[\mathscr{A},\Psi]t).
\end{equation*}
The definition of  $\Psi$ (\ref{psi2}) gives
\begin{equation}
e^{-\mathscr{A}t}F_0e^{\mathscr{A}t}=2b\mathscr{J}-2bt[\mathscr{A},\mathscr{J}]+\mathscr{A}^2+\frac{b-1}{N-1}\mathscr{A}^2\mathscr{J}-
\frac{b-1}{N-1}\mathscr{A}^2[\mathscr{A},\mathscr{J}]t.\label{sumando3}
\end{equation}
We now substitute (\ref{sumando I}), (\ref{sumando 2}) and (\ref{sumando3}) in the definition of $\xi$ (\ref{a verificar diagonal}) obtaining
\begin{align*}
\xi&= \frac{b-1}{N-1}(-\mathscr{A}\mathscr{J}\mathscr{A}-\mathscr{A}[\mathscr{A},\mathscr{J}]+\mathscr{A}^2\mathscr{J})
+2bt(\mathscr{A}-[\mathscr{A},\mathscr{J}]\\ &\quad\quad -\frac{b-1}{2b(N-1)}\mathscr{A}^2[\mathscr{A},\mathscr{J}])
+2bt^2\mathscr{D}+bI+2b\mathscr{J},
\end{align*}
(\ref{sumaprinc}) of lemma \ref{sumaciones} finally gives
$$
\xi=bI+2bt^2\mathscr{D} +2b\mathscr{J},
$$
which it is indeed a diagonal matrix.

This finishes the proof of Theorem 1.1.

\end{proof}

It remains to prove Lemma \ref{sumaciones}.

\begin{proof} (of Lemma \ref{sumaciones})

\textit{First step. Proof of (\ref{suma2s})}:

Using induction on $k$, we easily find that $[A^k,\mathscr{J}]=kA$, $k\geq1 $. This shows that
$\ad^{n}_A\mathscr{J}=0$, $n\geq 2$.
The definition of $\mathscr{A}$ (\ref{Acurs}) gives now (\ref{suma2s}).

\textit{Second step. Proof of (\ref{dfe0}), (\ref{dfe1}) and (\ref{dfe2})}:

First of all, the definition of $\Psi$ (\ref{psi2}) gives
\begin{equation}\label{ccc4}
[\mathscr{A},\Psi ]=\frac{b-1}{N-1}[\mathscr{A},\mathscr{J}].
\end{equation}
The definition of $\mathscr{A}$ (\ref{Acurs}) and (\ref{suma2s}) show that the matrices $\mathscr{A}$, $[\mathscr{A},\mathscr{J}]$ and
$e^{\mathscr{A}t}$ commute (they are linear combination of powers of $A$). We then have
\begin{align}
\label{cc1}
\ad ^n_{\mathscr{A}}(\Psi )&=0,\quad n\ge 2,\\
\label{cc2}
[\mathscr{A},\mathscr{J}]e^{\mathscr{A}t}&=e^{\mathscr{A}t}[\mathscr{A},\mathscr{J}].
\end{align}
Using (\ref{adformula}) and (\ref{cc1}) we get
$$
e^{\mathscr{A}t}\Psi e^{-\mathscr{A}t}=\sum_{n=0}^\infty \frac{t^n}{n!}\ad ^n_\mathscr{A}\Psi=\Psi+t[\mathscr{A},\Psi].
$$
The definition of $\Psi$ (\ref{psi2}) and $F_2$ (\ref{F2}) give now (\ref{dfe0}).

In a similar way, we have
$$
e^{-\mathscr{A}t}\Psi e^{\mathscr{A}t}=\Psi-t[\mathscr{A},\Psi].
$$
Using now the definition of $\mathscr{D}$ (\ref{D}), (\ref{ccc4}) and (\ref{cc2}) we have
\begin{align*}
e^{\mathscr{A}t}\mathscr{D} e^{-\mathscr{A}t}\Psi&=e^{\mathscr{A}t}\mathscr{D}(\Psi-t[\mathscr{A},\Psi]) e^{-\mathscr{A}t}\\
&=-\displaystyle \frac{b}{2}I-\frac{(b-1)t}{N-1}e^{\mathscr{A}t}\mathscr{D}e^{-\mathscr{A}t}[\mathscr{A},\mathscr{J}].
\end{align*}
This proves (\ref{dfe1}). The proof of (\ref{dfe2}) is similar.

\textit{Third step. Proof of (\ref{suma_2k+1})}:
 Using the definition of $\mathscr{A}$ in (\ref{Acurs}) and (\ref{suma2s}), we have to prove the following identity
$$
\sum_{j=0}^{\left[\frac{N}{2}\right]-1}\alpha_jA^{2j+1}\sum_{j=0}^{\left[\frac{N}{2}\right]-1}(2j+1)\alpha_jA^{2j+1}=
\frac{2b(N-1)}{1-b}\sum_{j=1}^{\left[\frac{N-1}{2}\right]}\alpha_j\frac{(2j)^j}{(2j+1)^{j-1}}A^{2j}.
$$
This is equivalent to prove
\begin{equation}\label{ult}
\sum_{j,s=0}^{\left[\frac{N}{2}\right]-1}\alpha_j\alpha_s(2s+1)A^{2(j+s+1)}=
\frac{2b(N-1)}{1-b}\sum_{j=1}^{\left[\frac{N-1}{2}\right]}\alpha_j\frac{(2j)^j}{(2j+1)^{j-1}}A^{2j}.
\end{equation}

Taking into account that $A$ is nilpotent,
\begin{equation}\label{alfas}
\alpha_s\alpha_j=\alpha_{j+s}\frac{(2j+1)^{j-1}(2s+1)^{s-1}}{(2(j+s)+1)^{j+s-1}}\left(\begin{array}{c}s+j\\s \end{array}\right),
\end{equation}
and writing $k=j+s$, we find
\begin{align}\label{igualdad}
&\sum_{j,s=0}^{\left[\frac{N}{2}\right]-1}\alpha_j\alpha_s(2s+1)A^{2(j+s+1)}\\\nonumber &\quad =
\sum_{k=0}^{\left[\frac{N-1}{2}\right]-1}\frac{\alpha_{k}}{(2k+1)^{k-1}}A^{2(k+1)}
\sum_{m=0}^{k}\left(\begin{array}{c}k\\m \end{array}\right)(2m+1)^m(2(k-m)+1)^{k-m-1}.
\end{align}

We now use Abel's binomial identity (see for instance \cite[p. 18]{Ro}): for $z,w\in\mathbb{C}$, $w\neq 0$,
\begin{equation}\label{abel}
\sum_{m=0}^k\left(\begin{array}{c}k\\m \end{array}\right)(m+z)^m(k-m+w)^{k-m-1}=w^{-1}(z+w+k)^k.
\end{equation}
Then
$$
\sum_{m=0}^{k}\left(\begin{array}{c}k\\m \end{array}\right)(2m+1)^m(2(k-m)+1)^{k-m-1}=2^k(1+k)^k,
$$
(\ref{igualdad}) now gives
$$
\sum_{j,s=0}^{\left[\frac{N}{2}\right]-1}\alpha_j\alpha_s(2s+1)A^{2(j+s+1)}
=\sum_{j=0}^{\left[\frac{N-1}{2}\right]-1}\alpha_{j}\frac{2^j(j+1)^j}{(2j+1)^{j-1}}A^{2(j+1)}.
$$

Using (\ref{alfas}) one has $$\alpha_j\alpha_1=\frac{\alpha_{j+1}(j+1)(2j+1)^{j-1}}{(2j+3)^j},\qquad \textrm{with}\ \alpha_1=\frac{1-b}{4b(N-1)}.$$

Thus,
$$
\sum_{j=0}^{\left[\frac{N-1}{2}\right]-1}\alpha_{j}\frac{2^j(j+1)^j}{(2j+1)^{j-1}}A^{2(j+1)}=\frac{2b(N-1)}{1-b}
\sum_{j=0}^{\left[\frac{N-1}{2}\right]-1}
\alpha_{j+1}\frac{2^{j+1}(j+1)^{j+1}}{(2j+3)^{j}}A^{2(j+1)}.
$$
This proves (\ref{ult}) and then (\ref{suma_2k+1}) as well.

\bigskip

\textit{Four step. Proof of (\ref{sumaprinc})}:
Taking into account (\ref{suma2s}) and (\ref{suma_2k+1}), we have to prove the following identity
$$
\sum_{j=1}^{\left[\frac{N}{2}\right]-1}2j\alpha_jA^{2j+1}=\sum_{j=0}^{\left[\frac{N}{2}\right]-1}\alpha_jA^{2j+1}
\sum_{j=1}^{\left[\frac{N-1}{2}\right]}\alpha_j\frac{(2j)^j}{(2j+1)^{j-1}}A^{2j}.
$$
This is equivalent to prove
\begin{equation}\label{sumatoria}
\sum_{j=1}^{\left[\frac{N}{2}\right]-1}2j\alpha_jA^{2j+1}=
\sum_{s=0}^{\left[\frac{N}{2}\right]-1}\sum_{j=1}^{\left[\frac{N-1}{2}\right]}\alpha_s\alpha_j\frac{(2j)^j}{(2j+1)^{j-1}}A^{2(s+j)+1}.
\end{equation}
Using (\ref{alfas}) once again,  we have
\begin{align*}
&\sum_{s=0}^{\left[\frac{N}{2}\right]-1}\sum_{j=1}^{\left[\frac{N-1}{2}\right]}\alpha_s\alpha_j\frac{(2j)^j}{(2j+1)^{j-1}}A^{2(s+j)+1}\\&\qquad =
\sum_{s=0}^{\left[\frac{N}{2}\right]-1}\sum_{j=1}^{\left[\frac{N-1}{2}\right]}\alpha_{s+j}
\frac{(2j)^j(2s+1)^{s-1}}{(2(j+s)+1)^{j+s-1}}\left(\begin{array}{c}s+j\\s \end{array}\right)A^{2(s+j)+1}.
\end{align*}
Writing $k=s+j$, and taking into account that $A$ is nilpotent of order $N$, we get for the right hand side of (\ref{sumatoria}) the expression
$$
\sum_{k=1}^{\left[\frac{N}{2}\right]-1}\frac{\alpha_{k}}{(2k+1)^{k-1}}A^{2k+1}\sum_{j=1}^{k}\left(\begin{array}{c}k\\j \end{array}\right)
(2j)^{j}(2(k-j)+1)^{k-j-1}.
$$
Writing $m=k-j$, one obtains
$$
\sum_{k=1}^{\left[\frac{N}{2}\right]-1}\frac{2^{k-1}\alpha_{k}}{(2k+1)^{k-1}}A^{2k+1}\sum_{m=0}^{k-1}\left(\begin{array}{c}k\\m \end{array}\right)
(k-m)^{k-m}(m+\frac{1}{2})^{m-1}.
$$
Abel's binomial identity (\ref{abel}) now gives
$$
\sum_{m=0}^{k-1}\left(\begin{array}{c}k\\m \end{array}\right)
(k-m)^{k-m}(m+\frac{1}{2})^{m-1}=\frac{1}{2^{k-1}}\left((2k+1)^k-(2k+1)^{k-1} \right).
$$
From where one can easily deduce (\ref{sumatoria}).

This proves (\ref{sumaprinc}).

The proof of the Lemma is now complete.

\end{proof}

\section{Rodrigues formula}

In this section we will prove Theorem \ref{th2} which provides a Rodrigues' formula for a sequence of orthogonal polynomials with respect to
the weight matrix $W$ for size $2\times 2$ (\ref{defW2}).

Write $R_n$ for the functions
\begin{equation}\label{defRn}
R_n=(-1)^ne^{- t^2}\left(\begin{array}{cc} \displaystyle{b^{-n} e^{(1-b) t^2}+ \frac{\vert a\vert ^2}{2}(n+2 t^2)}&a t\\
\bar a \left[2 t+e^{t^2} \sqrt{\pi}n\left( \Erf(\sqrt{b} t)-\Erf(t)\right)\right]&2
\end{array}
\right) ,
\end{equation}
where, as usual, $\Erf$ denotes the error function $\Erf(z)=\displaystyle{\frac{2}{\sqrt{\pi} }}\int_0^{z}e^{-x^2}dx$.

The Rodrigues' formula (\ref{defRodri}) can then be written as $P_n=R_n^{(n)}W^{-1}$.

First of all, we explain how one can use Theorem \ref{throdri}  to find these functions $R_n$, $n\ge 1$.

Indeed, Theorem \ref{th1} for size $2\times 2$ gives for the weight matrix $W$ the following symmetric second order
differential operator
\begin{equation}\label{difopepar}
D=\left( \frac{d}{dt}\right)^2\left(\begin{array}{cc} 1&a(b-1)t\\0&b
\end{array}
\right)+\left( \frac{d}{dt}\right)\left(\begin{array}{cc} -2bt&2ab\\0&-2bt
\end{array}
\right)+\left(\begin{array}{cc} 0&0\\0&2b
\end{array}\right) .
\end{equation}
Since the operator $D$ is symmetric with respect to $W$, the $n$-th  monic orthogonal polynomial $\hat P_n$ with respect to $W$ satisfies the differential equation $D(\hat P_n)=\Lambda_n\hat P_n$, where the eigenvalues $\Lambda_n$ are given by
$$
\Lambda_n=\left(\begin{array}{cc} -2bn&0\\0&-2b(n-1)
\end{array}
\right).
$$
Theorem \ref{throdri} associates the following
 second order differential equation (\ref{eqXn})  to the differential operator $D$
(\ref{difopepar}):
\begin{equation}\label{adifecu}
\left[R_n\left(\begin{array}{cc} 1&0\\a(b-1)t&b
\end{array}\right)\right]''-\left[R_n\left(\begin{array}{cc} -2bt&0\\a[b(2+n)-n]&-2bt
\end{array}\right)\right]'+R_n\Lambda_n = \Lambda_nR_n.
\end{equation}
Take now a solution $R_n$ of this differential equation and write $Y_n=R_n^{(n)}W^{-1}$.
Theorem \ref{throdri} guarantees that the function $Y_n$ satisfies the differential equation
$D(Y_n)=\Lambda _n Y_n$. Notice that $Y_n$ and $\hat P_n$ satisfy the same differential equation.
We have hence looked for a solution $R_ n$ of the differential equation (\ref{adifecu}) such that the matrix function
$R_n^{(n)}W^{-1}$ is also a matrix polynomial of degree $n$ with
nonsingular leading coefficient. This is the procedure we have used to find the functions $R_n$ given by (\ref{defRn}).

We now prove Theorem \ref{th2}, which establishes that actually the functions $R_n^{(n)}W^{-1}$ define
a sequence of orthogonal polynomials with respect to $W$.

\begin{proof} (of Theorem \ref{th2})

Using the Rodrigues' formula for Hermite polynomials, $H_n(t)=(-1)^n\left( e^{-t^2}\right)^{(n)}e^{t^2}$, \cite[Chapter 5]{Sz}, we have that
\begin{align*}
 \left( t^2e^{-t^2}\right)^{(n)} &=\frac{(-1)^n}{4}H_{n+2}(t)e^{-t^2}+\frac{(-1)^n}{2}H_{n}(t)e^{-t^2},\\
\left( te^{-t^2}\right)^{(n)} &=\frac{(-1)^n}{2}H_{n+1}(t)e^{-t^2},\\
\left( e^{-bt^2}\right)^{(n)}& =(-1)^n(\sqrt{b})^nH_{n}(\sqrt{b}t)e^{-bt^2},\\
(\Erf(t))^{(n)}& =(-1)^{n-1}\frac{2}{\sqrt{\pi}}H_{n-1}(t)e^{-t^2},\\
(\Erf(\sqrt{b}t))^{(n)}& =(-1)^{n-1}b^{n/2}\frac{2}{\sqrt{\pi}}H_{n-1}(\sqrt{b}t)e^{-bt^2}.
\end{align*}

These identities give, after straightforward computations using the three term recurrence relation for the Hermite
polynomials $tH_n=H_{n+1}/2+nH_{n-1}$ (see \cite[Chapter 5]{Sz}):
$$
 R_{n}^{(n)}(t)= \left(\begin{array}{cc}\displaystyle b^{-n/2}H_{n}(\sqrt{b}t)
 e^{-bt^2}+\frac{|a|^2}{2}tH_{n+1}(t)e^{-t^2} &\displaystyle \frac{a}{2}H_{n+1}(t)e^{-t^2}
 \\2\overline{a}\left(tH_{n}(t)e^{-t^2}-nb^{n/2}H_{n-1}(\sqrt{b}t)e^{-bt^2}\right)&2H_{n}(t)e^{-t^{2}}
\end{array}\right) .
$$
Taking into account that
$$
W^{-1}=\left(
\begin{array}{cc}
 e^{b t^2} & -a e^{b t^2} t \\
 -\overline{a} e^{b t^2} t & |a|^2 e^{b t^2} t^2+e^{t^2}
\end{array}
\right),
$$
we finally find the expression (\ref{PnitHn}) for $P_n$  given in Corollary \ref{cor1}:
$$
R_{n}^{(n)}(t)W^{-1}=\left(\begin{array}{cc}\displaystyle b^{-n/2}H_{n}(\sqrt{b}t)&-atb^{-n/2}H_{n}(\sqrt{b}t)+\frac{a}{2}H_{n+1}(t)\\
-2\overline{a}b^{n/2}nH_{n-1}(\sqrt{b}t)&2|a|^2b^{n/2}ntH_{n-1}(\sqrt{b}t)+2H_{n}(t)
\end{array}\right).
$$
This shows that $R_{n}^{(n)}(t)W^{-1}$ is a polynomial of degree $n$ (note that the entry $(1,2)$
of this matrix is
actually a polynomial of degree $n-1$) with nonsingular leading coefficient equal to $\Gamma _n$ (\ref{gamma}).

The orthogonality of $P_n$ and $t^kI$, $0\le k\le n-1$, with respect to $W$ follows taking into account that
$$
\int P_nW(t)t^{k}dt= \int R_n^{(n)}(t)t^kdt,
$$
and performing a careful integration by parts.

\end{proof}

\section{Three term recurrence relation}

\medskip
In order to find the recurrence coefficients (\ref{ttrrPn}) in the three term recurrence relation of Theorem \ref{the3}
we have followed the strategy of \cite{DG3} or \cite{DL}.

\begin{proof} (of Theorem \ref{the3})

We first compute the $L^2$ norm of the monic orthogonal polynomials $\hat P_n$ with respect to $W$.
Using the Rodrigues' formula (\ref{defRodri}), we have
$$
\Vert \hat P_n \Vert ^2=\int \hat P_n(t)W(t)t^ndt=\Gamma _n^{-1}\int R_n^{(n)}t^ndt,
$$
where $\Gamma _n$ are the leading coefficient of $P_n$ (\ref{gamma}). An integration by parts and the formulas
for $R_n$ and $\Gamma _n$ (see Theorem \ref{th2}) then give
\begin{equation}
\Vert \hat P_n \Vert ^2=\frac{\sqrt{\pi}n!}{2^{n}}\left(\begin{array}{cc} \displaystyle{\frac{\gamma_{n+1}}
{2b^{\frac{2n+1}{2}}}}&0\\0&\displaystyle\frac{2}{\gamma_n}
\end{array}
\right) .
\end{equation}
If we write
\begin{equation}\label{delta}
\Delta _n=\sqrt{\frac{2^{n}}{\sqrt{\pi}n!}}\left(\begin{array}{cc}
\displaystyle{\sqrt{\frac{2b^{\frac{2n+1}{2}}}{\gamma_{n+1}}}}&0\\0&\displaystyle{\sqrt{\frac{\gamma_n}{2}}},
\end{array}
\right),
\end{equation}
the polynomials
$$
\mathscr{P}_n=\Delta_n\hat P_n
$$
are then orthonormal with respect to $W$.

We now prove that the they satisfy the three term recurrence relation (\ref{ttrcon}).

This is just a matter of computation. Indeed, the coefficient $A_{n+1}$ in (\ref{ttrrPn}) is then
$$
A_{n+1}=\Delta_n \Delta _{n+1}^{-1}.
$$
The formula (\ref{delta}) for $\Delta_n$ gives now the formula for $A_n$ in (\ref{ttrrPn}).

On the other hand we have for the recurrence coefficient $B_n$ in (\ref{ttrcon}) the expression
$B_n=\Delta_n \hat B_n \Delta _{n}^{-1}$, where
$$
\hat B_n= \mbox{coeff. of $t^{n-1}$ in $\hat P_{n}-$ coeff. of $t^{n}$ in $\hat P_{n+1}$.}
$$
From (\ref{PnitHn}), we get that
 \begin{equation}\label{Bn monico}
\hat B_n=(b+(b-1)n)\left(\begin{array}{cc} 0&\displaystyle{\frac{a}{2b}}\\\displaystyle{\frac{2\bar{a}b^{n-\frac{1}{2}}}{\gamma_{n}\gamma_{n+1}}}&0
\end{array}\right),
\end{equation}
and the formula for $B_n$ in (\ref{ttrrPn}) follows easily.

\end{proof}

\bigskip

The three term recurrence relation for the polynomials $(P_n)_n$ (\ref{defRodri}) can now easily be computed.
Indeed, taking into account the expression for the leading coefficient $\Gamma _n$ (\ref{gamma}) of $P_n$ and $\Delta _n$ (\ref{delta}) of
$\mathscr{P}_n$, we find that $P_n= G_n \mathscr{P}_n$ where
\begin{equation}\label{gen}
G_n= \Gamma_n\Delta_n^{-1}=\sqrt{2^n\sqrt{\pi}n!}
\left(\begin{array}{cc}\displaystyle{\sqrt{\frac{\gamma_{n+1}}{2b^{\frac{2n+1}{2}}}}}&0\\0&\sqrt{2\gamma_n}
\end{array}\right).
\end{equation}

In particular, this gives $P_0=\displaystyle{\left(\begin{array}{cc} 1&0\\0&2
\end{array}
\right)}.$
If we write
\begin{equation}\label{recurr-normalizada}
tP_n(t)=\tilde A_{n+1}P_{n+1}(t)+\tilde B_nP_n(t)+\tilde C_nP_{n-1}(t),\quad n\geq 0,
\end{equation}
it follows from (\ref{ttrcon}) that $\tilde A_{n+1}=G_nA_{n+1}G^{-1}_{n+1}$, $\tilde B_{n}=G_nB_{n}G^{-1}_{n}$ and
$\tilde C_{n}=G_nA_{n}^*G^{-1}_{n-1}$. An easy computation gives now
\begin{align*}
\tilde A_n&=\frac{1}{2}\left(\begin{array}{cc} 1&0\\0&\displaystyle{\frac{\gamma_{n-1}}{\gamma_n}}
\end{array}
\right),\quad \tilde B_n=(-n+(n+1)b)\left(\begin{array}{cc} 0&\displaystyle{\frac{a}{2b\gamma_n}}
\\\displaystyle{\frac{2\bar{a}b^{n-\frac{1}{2}}}{\gamma_{n+1}}}&0
\end{array}
\right),\\
\tilde C_n&=n\left(\begin{array}{cc} \displaystyle{\frac{\gamma_{n+1}}{b\gamma_n}}&0\\0&1
\end{array}
\right).
\end{align*}

The $L^2$ norm of the polynomials $P_n$ follows easily from the formula (\ref{gen}) for the matrices $G_n$:
\begin{equation*}
\Vert P_n \Vert ^2=2^{n}\sqrt{\pi}n!\left(\begin{array}{cc} \displaystyle{\frac{\gamma_{n+1}}{2b^{n+\frac{1}{2}}}}&0\\0&2\gamma_n
\end{array}
\right) .
\end{equation*}

In a similar way, the three term recurrence relation for the monic orthogonal polynomials $(\hat P_n)_n$ can be deduced:
$$
t\hat P_n(t)=\hat P_{n+1}(t)+\hat B_n\hat P_n(t)+\hat C_n\hat P_{n-1}(t),\quad n\geq 0,
$$
where $\hat B_n$ is the one in (\ref{Bn monico}) and
$$
\hat C_n=\frac{n}{2b}\left(\begin{array}{cc} \displaystyle{\frac{\gamma_{n+1}}{\gamma_n}}&0\\0&\displaystyle{\frac{b\gamma_{n-1}}{\gamma_{n}}}
\end{array}
\right).
$$

\end{document}